\documentclass{amsart}

\usepackage[mathscr]{eucal}
\usepackage{amssymb}
\usepackage{amsthm}
\usepackage{amsmath}

\usepackage[all]{xy}
\CompileMatrices

\newdir{ >}{{}*!/-10pt/\dir{>}}

\newcommand{\uPt}[3]{#1}
\newcommand{\dPt}[3]{#2}
\newcommand{\tPt}[3]{#3}

\newcommand{\uPs}[6]{#1}
\newcommand{\dPs}[6]{#2}
\newcommand{\tPs}[6]{#3}
\newcommand{\qPs}[6]{#4}
\newcommand{\cPs}[6]{#5}
\newcommand{\sPs}[6]{#6}

\newtheorem{Thm}{Theorem}[section]
\newtheorem*{Thm*}{Theorem}

\newtheorem{Prop}[Thm]{Proposition}

\newtheorem{Lem}[Thm]{Lemma}
\newtheorem{Cor}[Thm]{Corollary}
\newtheorem*{Cor*}{Corollary}

\theoremstyle{definition}
\newtheorem{Def}[Thm]{Definition}
\newtheorem{Defs}[Thm]{Definitions}
\newtheorem{Not}[Thm]{Notation}

\newtheorem{Cons}[Thm]{Construction}

\theoremstyle{remark}
\newtheorem{Rem}[Thm]{Remark}
\newtheorem{Rems}[Thm]{Remarks}

\newtheorem{emptythm}[Thm]{}

\newcommand{\Disp}{\displaystyle}

\newcommand{\bbF}{\mathbb{F}}
\newcommand{\bbG}{\mathbb{G}}

\newcommand{\bbP}{\mathbb{P}}

\newcommand{\bbZ}{\mathbb{Z}}

\newcommand{\calO}{\mathcal{O}}

\newcommand{\isoto}{\buildrel \sim\over\to}
\newcommand{\too}{\mathop{\longrightarrow}}
\newcommand{\isotoo}{\buildrel \sim\over\too}
\newcommand{\onto}{\mathop{\twoheadrightarrow}}
\newcommand{\into}{\mathop{\rightarrowtail}}

\newcommand{\cat}[1]{\mathscr{#1}}

\newcommand{\cH}{\text{\rm\v H}}
\DeclareMathOperator{\id}{id}

\DeclareMathOperator{\Img}{Im}
\DeclareMathOperator{\Ker}{Ker}
\DeclareMathOperator{\Spec}{Spec}
\DeclareMathOperator{\Proj}{Proj}
\DeclareMathOperator{\Pic}{Pic} 
\newcommand{\Piccat}[1]{\Pic(\cat{#1})}
\newcommand{\Gm}{\bbG_{\rm m}}
\newcommand{\Gmcat}[1]{\Gm(\cat{#1})}
\newcommand{\cone}{{\operatorname{cone}}}

\newcommand{\Hom}{{\operatorname{Mor}}}
\newcommand{\homo}{\operatorname{hom}}
\newcommand{\iHom}{{\underline{\homo}}}
\newcommand{\Homcat}[1]{\Hom_{\cat #1}}

\newcommand{\inv}{^{-1}}

\newcommand{\op}{^{^{\text{op}}}}
\newcommand{\pullback}{\ar@{}[rd]|{{\displaystyle\lrcorner}}\ar@{} [rd]|{{\dot{}\,\ }}}
\newcommand{\smallbullet}{{\scriptscriptstyle\bullet}}

\newcommand{\equalby}[2]{\stackrel{#2}{#1}}
\newcommand{\equalbydef}{\equalby{=}{\text{def.}}}

\newcommand{\mm}{\mathfrak{m}}
\newcommand{\Mid}{\,\big|\,}

\newcommand{\restr}[1]{_{|_{\scriptstyle #1}}}
\newcommand{\oursetminus}{\smallsetminus}

\newcommand{\minifrac}[2]{\genfrac{}{}{0.4pt}{}{\scriptstyle #1} {\scriptstyle #2}}

\newcommand{\smallmatrice}[1]{\bigl(\begin{smallmatrix} #1 \end {smallmatrix}\bigr)}

\newcommand{\afortiori}{{\sl a fortiori}}
\newcommand{\eg}{{\sl e.g.}}

\newcommand{\ie}{{\sl i.e.}\ }

\newcommand{\loccit}{{\sl loc.\ cit.}}

\newcommand{\tcat}[1]{(\cat{#1},\otimes,\unit)}

\newcommand{\tenscat}{tensor triangulated category}
\newcommand{\Tenscat}[1]{\tenscat\ $\tcat{#1}$}

\newcommand{\tthsub}{thick $\otimes$-ideal}

\newcommand{\strongcl}{strongly closed}
\newcommand{\sctenscat}{\strongcl\ \tenscat}

\newcommand{\icsctenscat}{idempotent complete \strongcl\ \tenscat}

\newcommand{\iccat}[1]{\widetilde{\cat #1}}

\newcommand{\LU}[1]{\cat{L}(#1)}
\newcommand{\KU}[1]{\cat{K}(#1)}
\newcommand{\KZ}[1]{\cat{K}_{#1}}

\newcommand{\unit}{{1}}

\DeclareMathOperator{\Der}{D}
\DeclareMathOperator{\smallperf}{perf}
\newcommand{\Dperf}{\Der^{\smallperf}}
\DeclareMathOperator{\stab}{-\,stab}

\DeclareMathOperator{\Spc}{Spc}
\newcommand{\Spccat}[1]{\Spc(\cat #1)}

\DeclareMathOperator{\supp}{supp}

\DeclareMathOperator{\UU}{U}

\newcommand{\suppcat}[1]{\supp(\cat #1)}

\newcommand{\Ucap}{U_1\cap U_2}

\newcommand{\xytriangle}[7]{\xymatrix@C=#7em{#1\ar[r]^-{\Disp #4} &  #2 \ar[r]^-{\Disp #5}&#3\ar[r]^-{\Disp #6}&T #1}}

\newcommand{\xyTriangle}[8]{\xymatrix@C=#8em{#1\ar[r]^-{\Disp #5}&#2 \ar[r]^-{\Disp #6}&#3\ar[r]^-{\Disp #7}&#4}}

\newcommand{\degone}{{\scriptscriptstyle\cdot}}

\newcommand{\octahedron}[7]{\xymatrix@C=#6em @R=#7em{%
&&& \tPs #1 
\ar[rrdd]^-{\dPt #4} 
\ar[ld]|(.65){\dPt #3} 
\ar@{<-}[dd]|!{[rrdd];[ld]}{\hole} 
\\
\dPs #1 
\ar[rrru]^-{\cPs #1} 
\ar[rrdd]_-{\dPt #2} 
&& \uPt #3 
\ar[ll]_(.4){\tPt #3}|(.7){\degone} 
\ar[dd]_-{\tPt #5}|(.75){\degone} 
&
\\
&&& \uPs #1 
\ar[lllu]^(.6){\qPs #1}|!{[ld];[lu]}{\hole} 
\ar[uu] |!{[rr];[lu]}{\hole} _(.6){\sPs #1}
&& \uPt #4 
\ar[ll]_(.6){\tPt #4}|(.7){\degone} 
\ar[lllu]_-{\dPt #5} 
\\
&& \uPt #2 
\ar[ru]|(.55){\tPt #2}|(.75){\degone} 
\ar[rrru]_-{\ \uPt #5} 
}}

\begin{document}


\title{Gluing in tensor triangular Geometry}
\author{Paul Balmer}
\author{Giordano Favi}

\date\today

\address{Paul Balmer, Department Mathematics, ETH Zurich, 8092 Zurich, Switzerland}
\email{balmer@math.ethz.ch}
\urladdr{http://www.math.ethz.ch/$\sim$balmer}

\address{Giordano Favi, Department Mathematics, Universit\"at Basel, 4051 Basel, Switzerland}
\email{giordano.favi@unibas.ch}
\urladdr{http://www.math.unibas.ch/$\sim$giordano}

\begin{abstract}
We discuss gluing of objects and gluing of morphisms in tensor triangulated categories. We illustrate the results by producing, among other things, a Mayer-Vietoris exact sequence involving Picard groups.
\end{abstract}

\thanks{Research supported by Swiss National Science Foundation, grant~620-66065.}

\maketitle


\vskip-\baselineskip
\vskip-\baselineskip
\vskip-\baselineskip
\tableofcontents


\vskip-\baselineskip
\vskip-\baselineskip
\vskip-\baselineskip
\section*{Introduction}
\bigbreak

Tensor Triangular Geometry is the geometry of tensor triangulated categories. Heuristically, this contains at least Algebraic Geometry and the geometry of Modular Representation Theory but it also appears in many other areas of Mathematics, as recalled in the introduction of~\cite{b05}.

We will denote by $\cat{K}$ a triangulated category (with suspension $T:\cat{K}\isoto\cat{K}$) equipped with a tensor product, \ie an exact symmetric monoidal structure $\otimes: \cat{K}\times\cat{K}\too\cat{K}$, see more in Section~\ref{basics-sect}. Two key examples to keep in mind appear respectively in Algebraic Geometry, as $\cat{K}=\Dperf(X)$, the derived category of perfect complexes over a quasi-compact and quasi-separated scheme~$X$ (\eg\ a noetherian scheme), and in Modular Representation Theory, as $\cat{K}=kG\stab$, the stable category of finite dimensional representations modulo projective ones, for $G$ a finite group and $k$ a field of characteristic~$p>0$, typically dividing the order of the group.

In~\cite{b05}, the concept of \emph{spectrum} $\Spccat{K}$ of such categories is introduced. It is the universal topological space in which one can define \emph{supports} $\supp(a)\subset\Spccat{K}$ for objects $a\in\cat{K}$ in a reasonable way. In the above two examples, the spectrum $\Spccat{K}$ is respectively isomorphic to the scheme~$X$ itself and to the projective support variety~$\Proj\,\text{H}^{\smallbullet}(G,k)$.

\smallbreak

One fundamental construction of~\cite{b05} is the \textit{presheaf of triangulated categories}, $U\mapsto\KU{U}$, which associates to an open $U\subset\Spccat{K}$ a tensor triangulated category $\KU{U}$ defined as follows. Consider $Z=\Spccat{K}\oursetminus U$ the closed complement of~$U$ and consider the thick subcategory $\KZ{Z}\subset\cat{K}$ of those objects $a\in\cat{K}$ with $\supp(a)\subset Z$, \ie those objects which ought to disappear on~$U$. Then, the category
$$
\KU{U}:=\iccat{{}\cat{K}/\KZ{Z}}
$$
is defined as the idempotent completion of the Verdier quotient $\cat{K}/\KZ{Z}$. Localization $\cat{K}\onto\cat{K}/\KZ{Z}$ followed by idempotent completion $\cat{K}/\KZ{Z}\hookrightarrow\KU{U}$ yields a restriction functor $\rho_U:\cat{K}\to \KU{U}$. In the scheme example, it is an important theorem of Thomason~\cite{thtr} that for a quasi-compact open $U\subset X$ and for $\cat{K}=\Dperf(X)$, the above $\KU{U}$ is equivalent to $\Dperf(U)$. This is one reason for working with idempotent complete categories. Another reason is a key result of~\cite{b06} which says that if $\cat{K}$ is idempotent complete and if the support of an object of~$\cat{K}$ decomposes into two connected components then the object itself decomposes into two direct summands accordingly, see Theorem~\ref{Balmson} below.

\smallbreak

The present paper deals with the following type of questions. Suppose that $\Spccat{K}$ is covered by two quasi-compact open subsets $\Spccat{K}=U_1\cup U_2$ and consider the commutative diagram of restrictions\,:
\begin{equation}
\label{K12-eq}
\vcenter{\xymatrix{{\cat{K}}\ar[r]\ar[d]&{\KU{U_1}=:\cat{K}_1}\kern-3em\ar[d]
\\
\kern-3em{\cat{K}_2:=\KU{U_2}}\ar[r]&{\KU{U_1\cap U_2}=:\cat{K}_{12}}\,.\kern-4em
}}
\end{equation}
\textbf{Question\,:} Is the global category $\cat{K}$ obtained by ``gluing'' $\cat{K}_1$ and $\cat{K}_2$ over $\cat{K}_{12}$\,?

\smallbreak

This is a very natural question but it is known to be tricky, already in Algebraic Geometry. Indeed, it is easy to find non-zero morphisms $f:a\to b$ in $\cat{K}=\Dperf(X)$ such that $f\restr{U_1}=0$ and $f\restr{U_2}=0$ for an open covering $X=U_1\cup U_2$. Over $X=\bbP^1_{k}$ an example is the morphism $f:\calO(2)\to T(\calO)$ which is the third one in the exact triangle corresponding to the exact sequence $\calO\into \calO(1)\oplus\calO(1)\onto\calO(2)$\,; take for $U_1$ and $U_2$ two affine subsets. (For an exact sequence of vector bundles $E'\into E\onto E''$ over a scheme $X$, the corresponding morphism $f:E''\to T(E')$ is zero in $\Dperf(X)$ if and only if the sequence splits.) This example also shows that the phenomenon is not pathological but observable in very common situations.

\smallbreak

Still, the problem admits a nice solution, as explained in our main results\,:
\begin{Thm*}[\textbf{Mayer-Vietoris for morphisms}, see Thm.\,\ref{MV-mor-thm}]
In the above situation~\eqref{K12-eq}, given two objects $a,b\in\cat{K}$, there exists a long exact sequence\,:
$$
\cdots\,\Homcat{K_{12}}(Ta,b)\mathop{\to}^{\partial}\Homcat{K}(a,b)\to\Homcat{K_1}(a,b)\,\oplus\,\Homcat{K_2}(a,b)\to\Homcat{K_{12}}(a,b)\mathop{\to}^{\partial}\,\cdots
$$
The connecting homomorphism $\partial\,:\Homcat{K_{12}}(Ta,b)\too\Homcat{K}(a,b)$ is defined in Construction~\ref{partial-cons}. The other ones are the obvious restrictions and differences thereof.
\end{Thm*}

\begin{Thm*}[\textbf{Gluing of two objects}, see Thm.\,\ref{MV-obj-thm}]
In the above situation~\eqref{K12-eq}, given two objects $a_1\in\KU{U_1}$ and $a_2\in\KU{U_2}$ and an isomorphism $\sigma:a_1\isoto a_2$ in $\KU{U_1\cap U_2}$\,, there exists an object $a\in\cat{K}$ which becomes isomorphic to $a_i$ in $\KU{U_i}$ for $i=1,2$. Moreover, this gluing is unique up to (non-unique) isomorphism.
\end{Thm*}

We can extend the above result to three open subsets and three objects, at the cost of possibly losing unicity of the gluing\,:

\begin{Cor*}[\textbf{Gluing of three objects}, see Cor.\,\ref{MV-3obj-cor}]
Let $\Spccat{K}=U_1\cup U_2\cup U_3$ be a covering by quasi-compact open subsets. Consider three objects $a_i\in\KU{U_i}$ for $i=1,2,3$ and three isomorphisms $\sigma_{ij}:a_j\isoto a_i$ in $\KU{U_i\cap U_j}$ for $1\leq i<j\leq3$. Suppose that the cocycle relation $\sigma_{12}\circ\sigma_{23}=\sigma_{13}$ is satisfied in $\KU{U_1\cap U_2\cap U_3}$. Then there exists an object $a\in\cat{K}$, isomorphic to $a_i$ in $\KU{U_i}$ for $i=1,2,3$.
\end{Cor*}

In general, we do not know if this gluing is possible with more than three open subsets. Nevertheless, in Theorem~\ref{conn-glue-thm}, we give elementary conditions under which the gluing is possible for arbitrary coverings.

\smallbreak

Then, we apply the main results to obtain an exact sequence involving Picard groups. For us, the Picard group, $\Piccat{K}$, is the set of isomorphism classes of  invertible objects in $\cat{K}$, with the tensor product as multiplication. In Algebraic Geometry, $\Pic(\Dperf(X))$ is the usual Picard group of~$X$ up to possible shifts, see Prop.\,\ref{pic=pic-prop}. On the other hand, $\Pic(kG\stab)$ is nothing but the group of endo-trivial representations, which is one of the fundamental invariants of Modular Representation Theory. In the next statement, we denote by $\Gmcat{K}=\Homcat{K}(\unit,\unit)^\times$ the abelian group of automorphisms of the $\otimes$-unit object~$\unit\in\cat{K}$.

\begin{Thm*}[\textbf{Mayer-Vietoris for Picard groups}, see Thm.\,\ref{MV-pic-thm}]
Let $\Spccat{K}=U_1\cup U_2$ with $U_i$ quasi-compact. Then there is half a long exact sequence\,:
$$\xymatrix@C=2em@R=1em{&&\kern8em\cdots\ar[r]&\Homcat{K(\Ucap)}(T\unit,\unit)\ar[r]^-{1+ \partial}&
\\
\ar[r]^-{1+\partial}&{\Gmcat{K}}\ar[r]&{\Gm}(\KU{U_1})\oplus{\Gm}(\KU{U_2})\ar[r]
&{\Gm}(\KU{\Ucap})\ar[r]^-{\delta}&
\\
\ar[r]^-{\delta}&\Piccat{K}\ar[r]&\Pic(\KU{U_1})\oplus\Pic(\KU{U_2})\ar[r]& \Pic(\KU{\Ucap})\,.}$$
To the left, we have the Mayer-Vietoris long exact sequence, the homomorphism~$\partial$ is as before and the non-labelled morphisms are again the obvious restrictions and (multiplicative) differences of restrictions. The new homomorphism
$$
\delta\,:\Gm(\KU{\Ucap})\to\Piccat{K}
$$
assigns to a unit $\sigma\in\Gm(\KU{\Ucap})$ the invertible object obtained by gluing two copies of the objects $\unit\in\KU{U_1}$ and $\unit\in\KU{U_2}$ along $\sigma:\unit\isoto \unit$ on $\Ucap$.
\end{Thm*}

It would be very interesting to continue this sequence to the right, say, with Brauer groups of Azumaya algebras. Although this is still work in progress, the authors do not know yet whether such an extension is possible in general. Neither do we know what the Brauer group of $\cat{K}=kG\stab$ should be, for instance.

In fact, in Modular Representation Theory, applying the above results to $\cat{K}=kG\stab$ gives us a way to construct endo-trivial $kG$-modules from any \v Cech $\Gm$-cocycle over the projective support variety~$\Proj\,\text{H}^\smallbullet(G,k)$, as long as the involved covering has at most three open pieces. In particular, the map $\delta$ of the last result might be of interest to representation theorists and we do not know if it has been studied, even in special cases. Dave Benson and Jon Carlson suggested a possible link with the recent article~\cite{ca05}. This will be investigated in future work.

\smallbreak

Using the conditional gluing of more than three objects, we obtain the following result (Thm\,\ref{H1-thm}), which relates invertible modules over the spectrum $\Spccat{K}$ and invertible objects in~$\cat{K}$. See more comments in Remark~\ref{pic-rem}.

\begin{Thm*}
Suppose that $\Hom_{\KU{U}}(T\unit,\unit)=0$ for every quasi-compact open subsets~$U\subset\Spccat{K}$. Then, gluing induces an injective homomorphism from the first \v Cech cohomology of $\Spccat{K}$ with coefficients in~$\Gm$ into the Picard group of~$\cat{K}$
$$
\cH^1(\Spccat{K},\Gm)\hookrightarrow \Piccat{K}\,.
$$
\end{Thm*}

\smallbreak

We end the paper with the following formulation of Mayer-Vietoris\,:

\begin{Thm*}[\textbf{Excision}, see Thm.\,\ref{exc-thm}]
Let $Y\subset U\subset \Spccat{K}$. Assume that $Y$ is closed with quasi-compact complement and that $U$ is open and quasi-compact. Then the restriction functor $\cat{K}\to\KU{U}$ induces an equivalence between the subcategories of objects supported on $Y$, that is, $\KZ{Y}\isotoo\KU{U}_{Y}$.
\end{Thm*}


\bigbreak
\section{Basics about tensor triangulated categories and their geometry}
\label{basics-sect}
\medbreak


We survey the main concepts and results of~\cite{b05} and~\cite{b06}. Standard  notions about triangulated categories can be found in Verdier~\cite {verd} or Neeman~\cite{ne01}.

\begin{Defs}
\label{spc-def}
A \emph{tensor triangulated} category~$\tcat{K}$ is an essentially small triangulated category~$\cat{K}$ with a symmetric monoidal structure $\otimes:\cat {K}\times\cat{K}\too\cat{K}$, $(a,b)\mapsto a\otimes b$. We have in particular $a\otimes b\cong b\otimes a$ and $\unit \otimes a\cong a$ for the \emph{unit}~$\unit\in\cat{K}$. We assume moreover that the functors $a\otimes -$ and $-\otimes b$ are exact for every $a,b\in\cat {K}$ and that the usual diagram
$$
\xymatrix{T(a)\otimes T(b)\ar[r]^-{\cong}\ar[d]_-{\cong}&T(T(a)\otimes b)\ar[d]^-{\cong}\\T(a\otimes T(b))\ar[r]_-{\cong}&T^2(a\otimes b)}
$$
anticommutes. We use $T:\cat{K}\isoto\cat{K}$ to denote the translation (suspension).

A \emph{prime ideal} $\cat{P}\subsetneq\cat{K}$ is a proper  subcategory such that (1)-(4) below hold true\,:
\begin{enumerate}
\item $\cat{P}$ is a full triangulated subcategory, \ie $0\in\cat{P} $, $T(\cat{P})=\cat{P}$ and if  $a,b\in\cat{P}$ and if $a\to b\to c\to  T(a)$ is a distinguished triangle in~$\cat{K}$ then~$c\in\cat {P}$\,;
\smallbreak
\item $\cat{P}$ is thick, \ie if $a\oplus b\in\cat{P}$ then $a,b\in \cat{P}$\,;
\smallbreak
\item $\cat{P}$ is a $\otimes$-ideal, \ie if $a\in\cat{P}$ then
$a\otimes b\in\cat{P}$ for all $b\in\cat{K}$\,;
\smallbreak
\item $\cat{P}$ is \emph{prime}, \ie if $a\otimes b\in\cat{P}$ then $a \in\cat{P}$
or $b\in\cat{P}$.
\end{enumerate}

\goodbreak

A subcategory $\cat{J}\subset \cat{K}$ satisfying  (1), (2) and~(3) is a \emph{\tthsub}.

The \emph{spectrum} $\Spccat{K} $ is the set of primes $\cat{P}\subset\cat{K}$. The \emph{support} of  an object $a\in\cat{K}$ is defined as the subset $\supp(a)=\{\cat{P} \in\Spccat{K}\mid a\notin\cat{P}\}\subset\Spccat{K}$. The complements  $\UU(a)=\{\cat{P}\in\Spccat{K}\mid a\in\cat{P}\}$ of these supports  form a basis $\{\UU(a)\}_{a\in\cat{K}}$ of the so-called \emph {Zariski topology} on the spectrum.
\end{Defs}

\smallbreak

\begin{Thm}[\cite{b05} Thm.\,3.2]
\label{supp-thm}
Let $\cat{K}$ be a \tenscat. We have
\begin{enumerate}
\item[(i)] $\supp(0)=\varnothing$ and $\supp(\unit)=\Spccat{K}$.
\smallbreak
\item[(ii)] $\supp(a\oplus b)=\supp(a)\cup\supp(b)$.
\smallbreak
\item[(iii)] $\supp(Ta)=\supp(a)$.
\smallbreak
\item[(iv)] $\supp(a)\subset \supp(b)\cup\supp(c)$ for any distinguished
$a\to b\to c\to T(a)$.
\smallbreak
\item[(v)] $\supp(a\otimes b)=\supp(a)\cap \supp(b)$.
\end{enumerate}
Moreover, $(\Spccat{K},\supp)$ is universal for these properties.
\end{Thm}

\begin{Not}
\label{KZ-not}
Let $Y\subset \Spccat{K}$. We denote by $\KZ{Y}$ the full subcategory $\KZ{Y}:=\{a\in\cat{K}\mid\supp(a)\subset Y\}$ of  those objects \emph{supported on~$Y$}.
\end{Not}

\begin{Def}
\label{sc-def}
We call a \Tenscat{K} \emph{\strongcl} if there exists a bi-exact functor $\iHom:\cat{K}\op\times\cat{K}\too\cat{K}$ with natural  isomorphisms
\begin{equation}
\label{iHom-eq}
\Homcat{K}(a\otimes b,c)\cong\Homcat{K}(a,\iHom(b,c))
\end{equation}
and such that all objects are \emph{strongly dualizable}, \ie the  natural morphism
\begin{equation}
\label{sd-eq}
D(a)\otimes b\isoto\iHom(a,b)
\end{equation}
is an isomorphism for all $a,b\in\cat{K}$, where we denote by $D(a)$  the dual $D(a):=\iHom(a,\unit)$ of an object $a\in\cat{K}$. More  details can be found in~\cite[App.\,A]{hps}, for instance. It follows  from~\eqref{sd-eq} that $D^2(a)\cong a$ for all~$a\in\cat{K}$\,; see for instance~\cite[Thm.\,A.2.5\,(b)]{hps}.
\end{Def}

\begin{Prop}[\cite{b06} Cor.\,2.5]
\label{no-supp-prop}
Let $\cat{K}$ be a \sctenscat\ and let $a\in\cat{K}$ be an object. Then
$\supp(a)=\varnothing$ if and only if $a=0$.
\end{Prop}

\begin{Prop}[\cite{b06} Cor.\,2.8]
\label{Hom-prop}
Let $\cat{K}$ be a \sctenscat. Suppose that the supports of two  objects do not meet\,: $\supp(a)\cap\supp(b)=\varnothing$. Then there  is no non-trivial morphism between them\,: $\Homcat{K}(a,b)=0$.
\end{Prop}

\begin{Prop}
\label{iso-local-prop}
Let $\cat{K}$ be a \sctenscat. A morphism $f:a\to b$ in $\cat{K}$ is  an isomorphism if and only if it is an isomorphism in $\cat{K}/\cat{P} $ for all $\cat{P}\in\Spccat{K}$.
\end{Prop}

\begin{proof}
This easily follows from the fact that a morphism $f$ in a  triangulated category is an isomorphism if and only if $\cone(f)=0$.  Thus if $f$ is an isomorphism in $\cat{K}/\cat{P}$ we have that $\cone (f)\in\cat{P}$. If this is true for all $\cat{P}\in\Spccat{K}$ we  have that $\supp(\cone(f))=\varnothing$ which implies that $\cone(f)=0 $ by Proposition~\ref{no-supp-prop}.
\end{proof}

\begin{Thm}[\cite{b06} Thm.\,2.11]
\label{Balmson}
Let $\cat{K}$ be a \strongcl\ tensor triangulated category. Assume  that $\cat{K}$ is idempotent complete. Then, if the support of an  object $a\in\cat{K}$ can be decomposed as $\supp(a)=Y_1\cup Y_2$ for  disjoint closed subsets $Y_1$, $Y_2\subset\Spccat{K}$, with each open  complement $\Spccat{K}\oursetminus Y_i$ quasi-compact, then the  object itself can be decomposed as a direct sum $a\simeq a_1\oplus a_2 $ with $\supp(a_i)=Y_i$ for $i=1,2$.
\end{Thm}

\begin{Rem}
\label{ic-rem}
Recall that an additive category $\cat{K}$ is \emph{idempotent  complete} (or pseudo-abelian or karoubian) if all idempotents of all  objects split, that is, if $e\in\Homcat{K}(a,a)$ with $e^2=e$ then  the object $a$ decomposes as a direct sum $a\simeq a'\oplus a''$ on  which $e$ becomes $\smallmatrice{1&0\\0&0}$, that is, $a\simeq \Img(e) \oplus \Ker(e)$. One can always ``idempotent complete'' an additive  category~$\cat{K}\hookrightarrow\iccat{K}$ and $\iccat{K}$ inherits a unique triangulation from $\cat{K}$, see~\cite{basch}.
\end{Rem}

\begin{Defs}
\label{KU-def}
Let~$\cat{K}$ be an idempotent complete \sctenscat.

Let $U$ be a quasi-compact open subset of $\Spccat{K}$, and let us denote by $Z=\Spccat{K}\oursetminus U$ its closed complement. We denote by $\LU{U}=\cat{K}/\KZ{Z}$ the Verdier localization with respect to $\KZ{Z}$ (which can be realized by keeping the same objects as~$\cat{K}$ and by inverting all morphisms whose cone belongs to $\KZ{Z}$\,, by means of calculus of fractions). We denote by $\KU{U}=\iccat{{}\LU{U}}$ its idempotent completion. We have a fully faithful cofinal morphism $\LU{U}\hookrightarrow \KU{U}$ (\emph{cofinal} is sometimes called \emph{dense}, like in~\cite{thom}, and means that every object of the big category is a direct summand of an object of the small one).

For $U=\Spccat{K}$, by Proposition~\ref{no-supp-prop}, we have $\LU{U}=\cat{K}=\iccat{K}=\KU{U}$ since we assume~$\cat{K}$ idempotent complete. If $U\subset V$ we denote by $\rho_{V,U}:\LU{V}\to\LU{U}$ the localization functor and we also denote by
$$
\rho_{V,U}:\ \KU{V}\to\KU{U}
$$
the induced functor. When $V=\Spccat{K}$, we simply write $\rho_{U}:\cat{K}\to \KU{U}$ for $\rho_{V,U}$.

For two objects $a,b$ of $\cat{K}$ we denote by
$$
\Hom_U(a,b):=\Hom_{\LU{U}}\big(\rho_U(a),\rho_U(b)\big)=\Hom_{\KU{U}}\big(\rho_U(a),\rho_U(b)\big)
$$
the set of morphisms between $\rho_{U}(a)$ and $\rho_{U}(b)$ in $\LU{U}$ or equivalently in its idempotent completion $\KU{U}$\,; for simplicity, we might speak of ``morphisms between $a$ and $b$ in $\KU{U}$'', or simply of ``morphisms between $a$ and $b$ over~$U$''.
\end{Defs}

\begin{Prop}
\label{spc-KU-prop}
For $U\subset\Spccat{K}$ quasi-compact and open, the restriction functor $\rho_{U}:\cat{K}\to\KU{U}$ induces a homeomorphism $\Spc(\KU{U})\isoto U$, under which $\supp(\rho_{U}(a))=U\cap\supp(a)$, for any object $a\in\cat{K}$.
\end{Prop}

\begin{proof}
In fact, by~\cite{b05} Cor.\,3.14, $\Spc(\KU{U})=\Spc(\cat{K}/\KZ{Z})$ and by \loccit\ Prop.\,3.11, the localization functor induces a homeomorphism between $\Spc(\cat{K}/\KZ{Z})$ and the subspace $V:=\{\cat{P}\in\Spccat{K}\Mid \KZ{Z}\subset\cat{P}\}$ of~$\Spccat{K}$. So, it suffices to check that $V=U$. The last equality $\supp(\rho_{U}(a))=U\cap\supp(a)$ will then be a general fact about the functor $\Spc(-)$, see \loccit\ Prop.\,3.6.

Let $\cat{P}\in\Spccat{K}$. By the classification of thick $\otimes$-ideals, \loccit\ Thm.\,4.10, we have $\cat{P}\in V$, \ie $\KZ{Z}\subset \cat{P}$, if and only if $Z=\supp(\KZ{Z})\subset\suppcat{P}\equalbydef\cup_{a\in\cat{P}}\supp(a)$. By taking complements, this is equivalent to $\cap_{a\in\cat{P}}U(a)\subset U$, where $U(a)=\Spccat{K}\oursetminus\supp(a)=\{\cat{Q}\in\Spccat{K}\Mid a\in\cat{Q}\}$. Tautologically, $\cap_{a\in\cat{P}}U(a)=\{\cat{Q}\in\Spccat{K}\Mid\cat{P}\subset\cat{Q}\}$. The latter set is contained in~$U$ if and only if $\cat{P}\in U$\,: one direction is trivial and the other one uses that $Z$ is specialization closed, see \loccit\ Prop.\,2.9. So, $\cat{P}\in V$ if and only if $\cat{P}\in U$, as was left to check.
\end{proof}

\begin{Rem}
The above result cannot hold without assuming~$U$ quasi-compact since $\Spccat{K}$ is quasi-compact for any~$\cat{K}$. It is used above to insure $Z=\supp(\KZ{Z})$.
\end{Rem}

\bigbreak
\centerline{*\ *\ *}
\bigbreak

We end this Section with some general facts about triangulated categories.

\begin{Lem}
\label{four-tri}
Let $\cat{K}$ be a triangulated category. Then for every  distinguished triangle in which one object decomposes into two direct summands
$$\xytriangle{a}{b}{c_1\oplus c_2}{s}{\smallmatrice{\alpha\\\beta}} {\smallmatrice{\gamma&\,\delta}}{3}$$
there exist two objects, $d$ and $e$, and four distinguished triangles\,:
$$\xytriangle{d}{b}{c_1}{\alpha_0}{\alpha}{\alpha_2}{3}
\qquad
\xytriangle{a}{d}{c_2}{\delta_0}{\delta_1}{\delta}{3}$$
$$\xytriangle{e}{b}{c_2}{\beta_0}{\beta}{\beta_2}{3}
\qquad
\xytriangle{a}{e}{c_1}{\gamma_0}{\gamma_1}{\gamma}{3}$$
where $\alpha_2=T\delta_0\,\gamma$\,, $\delta_1=\beta\,\alpha_0$\,, $ \beta_2=T\gamma_0\,\delta$ and $\gamma_1=\alpha\,\beta_0$. Moreover, the morphism~$s$ factors as $s=\alpha_0\,\delta_0=\beta_0\,\gamma_0$.

In particular, we have $\cone(\alpha)\simeq\cone(\delta)$ and $\cone(\beta)\simeq\cone(\gamma)$.
\end{Lem}

\begin{proof}
We will prove the existence of the first two triangles, the other two  are obtained symmetrically ($c_1\oplus c_2\simeq c_2\oplus c_1$). The triangles are obtained by applying the Octahedron Axiom to the  equality $(1\ 0)\smallmatrice{\alpha\\\beta}=\alpha$ as displayed below:

$$\octahedron{b{c_1\oplus c_2}{c_1}{\smallmatrice{\alpha\\\beta}}{(1\  0)}{\alpha}}{{Ta}{(\gamma\ \delta)}{-Ts}}{{Tc_2}0{\smallmatrice{0\\  1}}}{{Td}{\alpha_2}{-T\alpha_0}}{{T\delta_0}{-T\delta_1}{T\delta}} {3.2}{3.2}$$
\end{proof}

\begin{Def}
\label{bicar-def}
We say that a commutative square as follows is a \textit{weak push-out}
$$
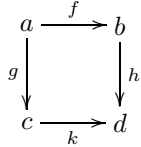
\xymatrix{a\ar[r]^-{f}\ar[d]_-{g}&b\ar[d]^-{h}\\c\ar[r]_-{k}&d}
$$
if there exists a distinguished triangle $\xytriangle{a}{b\oplus c}{d}{\smallmatrice{f\\g}}{\smallmatrice{-h&k}}{l}{2.5}$ for some morphism $l:d\to T(a)$. This is justified since $(d,h,k)$ satisfies the universal property of the push-out of $f$ and $g$ but without unicity of the factorization. Since such a square is then also a \textit{weak pull-back}, we call it \textit{weakly bicartesian}.
\end{Def}


\bigbreak
\section{Mayer-Vietoris for morphisms}
\bigbreak


%
\begin{Def}
\label{MV-sit-def}
Let $\cat{K}$ be an idempotent complete \sctenscat. We say that we are in \emph{a Mayer-Vietoris situation} if the spectrum of~$\cat{K}$ is covered by two quasi-compact open subsets $\Spccat{K}=U_1\cup U_2$. 

We shall denote by $Z_i=\Spccat{K}\oursetminus U_i$ the closed complements for $i=1,2$. Recall the important Definitions~\ref{KU-def}. We will use the simplified notation $\rho_i=\rho_{U_i}$ and $ \rho_{ij}=\rho_{\,U_i\,,\,U_i\cap U_j}$ for the restriction functors. We have the commutative diagram\,:
$$
\xymatrix@C=2em@R=2em{{\cat{K}}\ar@{->>}[rrd]_-{\Disp\rho_1}\ar@{->>}[rdd]^(.55){\Disp\rho_2}\ar[rrr]^-{\Disp\rho_1}\ar[ddd]_-{\Disp\rho_2}
&&&{\KU{U_1}}\ar[ddd]^-{\Disp\rho_{12}}
\\
&&{\LU{U_1}}\ar@{->>}[d]^-{\Disp\rho_{12}}\ar@{^(->}[ru]
\\
&{\LU{U_2}}\ar@{->>}[r]_-{\Disp\rho_{21}}\ar@{^(->}[ld]
&{\LU{U_1\cap U_2}}\ar@{_(->}[rd]
\\
{\KU{U_2}}\ar[rrr]_-{\Disp\rho_{21}}&&&{\KU{U_1\cap U_2}}
}
$$
where $\onto$ denotes a Verdier localization and $\hookrightarrow$ a fully faithful cofinal embedding.
\end{Def}

\begin{Rem}
\label{more-rem}
We do not really need to have $U_1$ and $U_2$ open and it would be enough to assume that they are arbitrary intersections of quasi-compact open subsets. Indeed, the key result taken from~\cite{b06}, Theorem~\ref{Balmson}, holds in this generality. Therefore, everything below holds in similar generality. We stick to open pieces because this is closer to the reader's understanding of a Mayer-Vietoris framework.
\end{Rem}

\begin{Def}
\label{U-iso-def}
Let $U\subset\Spccat{K}$ be a (quasi-compact) open with closed complement~$Z$. A morphism $s:a\to b$ in $\cat{K}$ is called a \emph{$U$-isomorphism} if it is an isomorphism in $\LU{U}$, or equivalently in~$\KU{U}$. This is also equivalent to saying that $\cone(s)$ belongs to $\KZ{Z}$ which also reads $\supp(\cone(s))\cap U= \varnothing$.
\end{Def}

\begin{Lem}
\label{U-iso-lem}
Consider a weakly bicartesian square in~$\cat{K}$ (Def.\,\ref{bicar-def})\,:
$$
\xymatrix{a\ar[r]^-{f}\ar[d]_-{g}&b\ar[d]^-{h}\\c\ar[r]_-{k}&d\,.\kern-0.5em}
$$
Let $U\subset\Spccat{K}$ be open. Then $f$ is a $U$-isomorphism if and only if $k$ is.
\end{Lem}

\begin{proof}
There exists a distinguished triangle $\xyTriangle{T\inv d}{a}{b\oplus c}{d\,.}{}{\smallmatrice{f\\g}}{\smallmatrice{-h&k}}{2.5}$ By Lemma~\ref{four-tri}, $\cone(f)\simeq\cone(k)$ and the result follows.
\end{proof}

\begin{Rem}
\label{X-iso-rem}
In a Mayer-Vietoris situation, a morphism which is both a $U_1$- and a $U_2$-isomorphism must be an isomorphism since the support of its cone is empty.
\end{Rem}

\begin{Lem}
\label{factorize-lem}
In a Mayer-Vietoris situation, suppose that $s:a\to b$ is a  $\Ucap$-isomorphism. Then $s$ may be factored as $s=s_1\, s_2$ where  $s_i$ is a $U_i$-isomorphism.
\end{Lem}

\begin{proof}
By hypothesis we have that $\cone(s)\in\KZ{Z_1\cup Z_2}$. Thus by  Theorem~\ref{Balmson} $\cone(s)$ may be written as $\cone(s)\simeq c_1 \oplus c_2$ where $c_i\in\KZ{Z_i}$. Now use Lemma \ref{four-tri}  which tells that $s=\alpha_0\,\delta_0$ and that $\cone(\alpha_0) \simeq c_1$ and $\cone(\delta_0)\simeq c_2$.
\end{proof}

\begin{Rem}
One can actually prove that the above factorization is essentially unique but we shall not use this fact below.
\end{Rem}

\begin{Lem}
\label{bicar-lem}
In a Mayer-Vietoris situation, consider a commutative diagram\,:
$$
\xymatrix{a\ar[r]^-{s_1}\ar[d]_-{t_2}&b\ar[d]^-{s_2}\\c\ar[r]_-{t_1}&d\,.}
$$
Assume that $s_i$ and $t_i$ are $U_i$-isomorphisms for $i=1,2$. Then the square is weakly bicartesian.
\end{Lem}

\begin{proof}
Consider the weak push-out $(e,u_1,u_2)$ of $s_1$ and $t_2$ and the morphism $v:e\to d$ induced by $s_2$ and $t_1$\,:
$$
\xymatrix{a\ar[r]^-{s_1}\ar[d]_-{t_2}&b\ar[d]_-{u_2}\ar@/^1em/[rd]^-{s_2}\\c\ar[r]^-{u_1}\ar@/_1em/[rr]_-{t_1}
&e\ar@{-->}[r]^-{v}
&d\,.}
$$
By Lemma~\ref{U-iso-lem}, $u_i$ is a $U_i$-isomorphism for $i=1,2$. By 2-out-of-3, $v$ is both a $U_1$- and a $U_2$-isomorphism, hence an isomorphism (see Rem.\,\ref{X-iso-rem}).
\end{proof}

\begin{Cons}
\label{partial-cons}
Consider a Mayer-Vietoris situation (Def.\,\ref{MV-sit-def}) and two objects $a,b\in\cat{K}$. We define a homomorphism
\begin{align*}
\partial\,:\ \Hom_{U_1\cap U_2}(a,b)&\too\Homcat{K}\big(a,T(b)\big)\\
g\quad&\longmapsto\quad\partial(g)
\end{align*}
as follows. Let $f\,s\inv=\big(\xymatrix@C=2em{a& x\ar[l]_-{s}\ar[r]^-{f} & b}\big)$  be a fraction representing a morphism $g\in\Hom_{\Ucap}(a,b)$. The cone of the $\Ucap$-isomorphism $s$ belongs to $\KZ{Z_1\cup Z_2}=\KZ{Z_1}\oplus\KZ{Z_2}$, see Theorem~\ref{Balmson}. So, we can chose a distinguished triangle
$$
\xytriangle{x}{a}{c_1\oplus c_2}{s}{\smallmatrice{\alpha\\ \beta}} {\smallmatrice{\gamma& \delta}}{3}
$$
where $c_i\in\KZ{Z_i}$ for $i=1,2$. Note that $\gamma\alpha+\delta\beta=0$. Define now
$$
\partial(g):=T(f)\circ \gamma\circ\alpha =-T(f)\circ\delta\circ\beta\,.
$$
This is a morphism in $\Homcat{K}(a,b)$, independant of the choices, see Theorem~\ref{MV-mor-thm}.
\end{Cons}

\begin{Rem}
\label{partial-rem}
Since $T$ is an equivalence, the above construction also induces\,:
$$
\xymatrix@C=2em{\Hom_{\Ucap}\big(T(a),b\big)\ar[r]^-\partial
&\Homcat{K}\big(T(a),T(b)\big)\ar[r]_-{\cong}^-{T\inv}
&\Homcat{K}(a,b)
}
$$
and we also denote this homomorphism by~$\partial$, since no confusion should follow. Explicitly, for a morphism $g=\big(\xymatrix@C=2em{T(a)& x\ar[l]_-{s}\ar[r]^-{f} & b}\big)$ in $\Hom_{\Ucap}(T(a),b)$, choose any distinguished triangle $\xyTriangle{a}{c_1\oplus c_2}{x}{T(a)}{\smallmatrice{\alpha\\ \beta}} {\smallmatrice{\gamma& \delta}}{s}{2}$ with $c_i\in\KZ{Z_i}$ for $i=1,2$\,; then we have $\partial(g)=f\circ\gamma\circ\alpha=-f\circ\delta\circ\beta\in\Homcat{K}(a,b)$.
\end{Rem}

\begin{Thm}[Mayer-Vietoris for morphisms]
\label{MV-mor-thm}
In a Mayer-Vietoris situation (Def.\,\ref{MV-sit-def}), for each pair of objects $a,b\in\cat{K} $, the above connecting homomorphism $\partial$ is well-defined and there is a long exact sequence of abelian groups
$$\xymatrix@C=1.8em @R=1em{&&\kern10em{\cdots}\ar[r]
&\Hom_{U_1\cap U_2}(Ta,b)\ar[r]^-{\partial} &
\\
\ar[r]^-{\partial}
&\Homcat{K}(a,b)\ar[r]
&\Hom_{U_1}(a,b)\oplus\Hom_{U_2}(a,b)\ar[r]
&\Hom_{U_1\cap U_2}(a,b)\ar[r]^-{\partial}&
\\
\ar[r]^-{\partial}
&\Homcat{K}(a,Tb) \ar[r]
&{\cdots}\kern10em
}$$
where the non-labelled homomorphisms are the obvious restrictions and differences of restrictions.
\end{Thm}

\begin{proof}
First, we have to check that the definition of $\partial(\minifrac{f}{s})$ given in~\ref{partial-cons} does not depend on the choice of the objects $c_i\in\KZ{Z_i}$ such that $\cone(s)\simeq c_1\oplus c_2$. This is easy, for other $c_i$ must be isomorphic to the chosen ones\,: $c_1\oplus c_2\simeq c_1'\oplus c_2'$ and $\supp(c_i)\subset Z_i$ forces $c_i\simeq c_i'$ for $i=1,2$, by Prop.\,\ref{Hom-prop}. The isomorphism $c_1\simeq c_1'$ disappears in the composition $\gamma\circ\alpha$ and \afortiori\ in $\partial(\minifrac{f}{s})$.

Then, we have to check that $\partial(\minifrac{f}{s})$ only depends on the class of the fraction $\minifrac{f}{s}$. To see this, let $t:y\to x$ be a $\Ucap$-isomorphism and let $\minifrac{ft}{st}=\big(\xymatrix@C=2em{a& y\ar[l]_- {st}\ar[r]^-{ft} & b}\big)$ be the amplified fraction. We have to show that $\partial(\minifrac{ft} {st})=\partial(\minifrac{f}{s})$. The morphism $st$ fits in a distinguished triangle.
$$
\xytriangle{y}{a}{d_1\oplus d_2}{st}{\smallmatrice{\alpha'\\ \beta'}} {\scriptstyle(\gamma'\, \delta')}{3}
$$
where $d_i\in\KZ{Z_i}$. Comparing the triangles for $s$ and for $st$ yields the diagram
$$\xymatrix@C=3em{y\ar[r]^-{\Disp st}\ar[d]_{\Disp t} & a\ar@{=}[d] \ar[r]^- {\smallmatrice{\alpha'\\ \beta'}}&d_1\oplus d_2\ar@{-->}[d]_{\Disp\exists\,\epsilon}\ar[r]^-{\scriptstyle(\gamma'\, \delta')}&T y\ar[d]^{\Disp Tt}\\
            x\ar[r]_-{\Disp s} & a \ar[r]_-{\smallmatrice{\alpha\\  \beta}}&c_1\oplus c_2\ar[r]_-{\scriptstyle(\gamma\, \delta)}&T x}$$
for some morphism $\epsilon$. But since $\supp(d_i)\cap\supp(c_j)=\varnothing$ for $i \neq j$, by Proposition~\ref{Hom-prop}, we have that $\epsilon=\smallmatrice {\epsilon_1 & 0 \\ 0 & \epsilon_2}$. Now compute
$$\partial\left(\minifrac{ft}{st}\right)=T(ft)\,\gamma'\,\alpha'=Tf\, Tt\, \gamma'\,\alpha'=Tf\,\gamma\,\epsilon_1\,\alpha'=Tf\,\gamma\,\alpha= \partial\left(\minifrac{f}{s}\right).$$
This proves that $\partial$ is well-defined. The fact that $\partial$  does not depend on the amplification of the fraction shows also that  in order to prove that $\partial$ is a group homomorphism it suffices  to see that $\partial(\minifrac{f+g}{s})=\partial(\minifrac{f}{s})+ \partial(\minifrac{g}{s})$, which is immediate.

\smallbreak

We now prove that the sequence is exact. It is easy to see that all consecutive compositions are zero. (Recall the notation for the restriction functors $\rho_i$ and $\rho_{ij}$ from Definition~\ref{MV-sit-def}.) For instance, $\rho_i(\partial(\minifrac{f}{s}))=0$ because $\partial(\minifrac{f}{s})$ factors via $c_i\in\KZ{Z_i}$ which becomes zero over~$U_i$\,. To see that $\partial\circ\rho_{ij}=0$, we check that $\partial(\minifrac{f}{s})=0$ if $s$ is a $U_1$-isomorphism for instance. But in this case, $c_2=0$ and $\partial(\minifrac{f}{s})$ factors via $c_2$\,.

\medbreak
\noindent\textit{Exactness at $\Hom_{U_1}(a,b)\oplus\Hom_{U_2}(a,b)$\,}:
Let $(f_1,f_2)\in \Hom_{U_1}(a,b)\oplus\Hom_{U_2}(a,b)$ such that $\rho_{12}(f_1)=\rho_{21}(f_2)$. Write $f_i=\big(\xymatrix@C=2em{a& x_i\ar[l]_-{s_i}\ar[r]^-{g_i} & b}\big)$.  Then there exist an object $x$ and $\Ucap$-isomorphisms $t_i:x\to x_i $ such that the diagram
\begin{equation}
\label{mor-1-eq}
\vcenter{
\xymatrix@C=2em@R=2em{ & x_1\ar[dl]_-{s_1}\ar[dr]^-{g_1} & \\ a&x\ar [u]_-{t_2}\ar[d]^-{t_1} & b\\ & x_2\ar[ul]^-{s_2}\ar[ur]_-{g_2}&}
}
\end{equation}
is commutative. By Lemma \ref{factorize-lem} we know that every $\Ucap$-isomorphism factors as a $U_1$-isomorphism followed by a $U_2$-isomorphism (and viceversa) so that we may choose $t_2$ to be a $U_2$-isomorphism, up to possibly amplifying the fraction $f_1$ without changing it. Similarly, we can assume $t_1$ is a $U_1$-isomorphism. By Lemma~\ref{bicar-lem}, the left ``square'' of~\eqref{mor-1-eq} is weakly bicartesian. Therefore (weak push-out), $g_1$ and $g_2$ induce a morphism $f:a\to b$ such that $f\circ s_i=g_i$ for $i=1,2$. Hence $f=g_i\, s_i^{-1}=f_i$ on $U_i$ as wanted.

\medbreak

\noindent\textit{Exactness at $\Homcat{K}(a,b)$\,}: Let $f\in\Homcat{K}(a,b)$ be such that $\rho_i(f)=0$ in $\LU{U_i}=\cat{K}/\KZ{Z_i}$ for $i=1,2$. This means that $f$ factors through objects $c_i\in\KZ{Z_i}$ as follows\,:
$$\xymatrix@C=2em@R=2em{ & c_1\ar[dr]^{f_1} & \\ a\ar[rr]^f\ar[ur]^{\alpha}\ar[dr]_{\beta}& & b\,.\kern-0.5em\\ & c_2\ar[ur]_{f_2}&}$$
Take now $x$ the weak push-out of $\alpha$ and $\beta$. By construction of the weak push-out (Def.\,\ref{bicar-def}), we have a distinguished triangle as in the first line of the diagram below. Since $\smallmatrice{f_1&\ -f_2}\cdot\smallmatrice{\alpha\\\beta}=f-f=0$, there exists a morphism $h:x\to b$ as follows\,:
$$\xymatrix{a\ar[r]^-{\smallmatrice{\alpha\\\beta}}& c_1\oplus c_2\ar[r]^-{\smallmatrice{\gamma&\delta}}\ar[dr]_{\smallmatrice{f_1&\,-f_2}}& x\ar@{-->}[d]^-{h}\ar[r]^{s} & Ta\\ & & b\,.\kern-0.5em}$$
We obtain a morphism $\minifrac{h}{s}=\big(\xymatrix{Ta & x\ar[l]_-{s}\ar[r]^-{h} & b}\big)\in\Hom_{\Ucap}(Ta,b)$. By Construction~\ref{partial-cons} and Remark~\ref{partial-rem}, we have $\partial(\minifrac{h}{s})=h\,\gamma\,\alpha=f_1\,\alpha=f$.

\medbreak

\noindent\textit{Exactness at $\Hom_{\Ucap}(a,b)$\,}: Let $\minifrac{f} {s}=\big(\xymatrix@C=2em{a& x\ar[l]_-{s}\ar[r]^-{f} & b}\big)$ be a morphism over $\Ucap$ such that $\partial(\minifrac{f}{s})=0$. As in Construction~\ref{partial-cons}, choose a distinguished triangle
$$
\xytriangle{x}{a}{c_1\oplus c_2}{s}{\smallmatrice{\alpha\\ \beta}} {\smallmatrice{\gamma&\delta}}{3}
$$
with $c_i\in\KZ{Z_i}$\,. The assumption $\partial(\minifrac{f}{s})=0$ reads $Tf\,\gamma\,\alpha=0$. Now apply Lemma~\ref{four-tri} to the above triangle to produce objects $d,e\in\cat{K}$ and morphisms $\alpha_0$, $\beta_0$, $\gamma_0$ and $\delta_0$ satisfying all the conclusions of Lemma~\ref{four-tri}, which the reader is encouraged to have at hand.

\textit{Claim}\,: There exists a distinguished triangle of the form
\begin{equation}
\label{tri-eq}
\xyTriangle{b}{Ta}{Td\oplus Te}{Tb\,.}{\gamma\,\alpha}{\smallmatrice{T \delta_0\\ -T\gamma_0}}{}{3}
\end{equation}
Indeed, the composition $\alpha_2\,\gamma_1=T\delta_0\,\gamma\,\gamma_1=0$ yields an Octahedron\,:
$$\octahedron{{e}{c_1}{Td}{\gamma_1}{\alpha_2}{0}}
{{Ta}{\gamma}{-T\gamma_0}} 
{{Tb} {-T\alpha_0}{-T\alpha}} 
{{Td\oplus Te} {\smallmatrice{1\\ 0}} {(0\ 1)}} 
{{\smallmatrice{\varphi\\ -T\gamma_0}}{\zeta}{-T(\gamma\alpha)}}
{3.2}{3.2}$$
for some morphisms $\varphi$ and $\zeta$. Note in particular the distinguished triangle
$$
\xyTriangle{b}{Ta}{Td\oplus Te}{Tb\,.}{\gamma\,\alpha}{\smallmatrice {\varphi\\ -T\gamma_0}}{\zeta}{3}
$$
To obtain the triangle~\eqref{tri-eq}, observe that $(T \delta_0-\varphi)\gamma=\alpha_2-\alpha_2=0$. By the distinguished triangle over~$\gamma$, there exists a morphism $h:Te\to Td$ such that $T \delta_0-\varphi=h\, T\gamma_0$. Using this equality we get an isomorphism of triangles
$$\xymatrix@C=3em{b\ar[r]^-{\gamma\,\alpha}& Ta\ar[r]^-{\smallmatrice {\varphi\\ -T\gamma_0}} & Td\oplus Te\ar[r]^-{\zeta}\ar[d]^{\simeq}_ {\smallmatrice{1 & -h\\ 0 & 1}} & Tb \\
b\ar@{=}[u]\ar[r]_-{\gamma\,\alpha} & Ta\ar@{=}[u]\ar[r]_- {\smallmatrice{T\delta_0\\ -T\gamma_0}} & Td\oplus Te\ar[r]_-{\zeta'} & Tb\ar@{=}[u]}$$
for $\zeta':=\zeta\cdot\smallmatrice{1 & h\\ 0 & 1}$. So, the lower triangle is distinguished. Hence the Claim.

Using this triangle and the assumption $Tf\circ\gamma\,\alpha=0$ yields a factorisation of $Tf$ as follows\,:
$$
\xymatrix{a\ar[r]^-{\gamma\,\alpha}
& Tx\ar[dr]_{Tf}\ar[r]^- {\smallmatrice{T\delta_0 \\ -T\gamma_0}}
& Td\oplus Te\ar@{-->}[d]^- {\smallmatrice{Tf_1&\,Tf_2}}\ar[r]
& Ta 
\\
&& Tb &
}
$$
for some morphisms $f_1:d\to b$ and $f_2:e\to b$. This reads $f=f_1\,\delta_0-f_2\, \gamma_0$. Using the triangles of Lemma~\ref{four-tri}, it is easy to see that $\alpha_0$, $\gamma_0$ are $U_1$-isomorphisms and that $\beta_0$ and $\delta_0$ are $U_2$-isomorphisms. Consider now the morphisms $\minifrac{f_1}{\alpha_0}=\big(\xymatrix@C=2em{a& d\ar[l]_-{\alpha_0}\ar[r]^-{f_1} & b}\big)$ and $ \minifrac{f_2}{\beta_0}=\big(\xymatrix@C=2em{a& e\ar[l]_-{\beta_0}\ar[r]^- {f_2} & b}\big)$ in $\LU{U_1}$ and $\LU{U_2}$ respectively. When restricted to $\LU{\Ucap}$ they clearly satisfy $\minifrac{f_1} {\alpha_0}-\minifrac{f_2}{\beta_0}=\minifrac{f_1\delta_0}{\alpha_0\delta_0}-\minifrac{f_2\gamma_0}{\beta_0\gamma_0}=\minifrac{f}{s}$. The last equation uses the relation $s=\alpha_0\delta_0=\beta_0\gamma_0$ from Lemma~\ref{four-tri}.
\end{proof}


\bigbreak
\section{Gluing objects}
\medbreak


It is convenient to fix the following standard terminology.

\begin{Def}
\label{glue-def}
Let $\Spccat{K}=U_1\cup \ldots\cup U_n$ be a covering by quasi-compact open subsets. Consider objects $a_i\in\KU{U_i}$ and isomorphisms $\sigma_{ji}:a_i\isoto a_j$ on $U_i\cap U_j$ such that $\sigma_{ki}=\sigma_{kj}\,\sigma_{ji}$ on $U_i\cap U_j\cap U_k$ for $1\leq i,j,k \leq n$. A \emph{gluing of the objects $a_i$ along the isomorphisms $\sigma_{ij}$} is an object $a\in\cat{K}$ and $n$ isomorphisms $f_i:a\isoto a_i$ on $U_i$ such that $\sigma_{ji}\,f_i=f_j$ on $U_i\cap U_j$ for all $1\leq i,j\leq n$. An \emph{isomorphism} of gluings $f:(a,f_1,\ldots,f_n)\isoto (a',f_1',\ldots,f_n')$ is an isomorphism $f:a\isoto a'$ in $\cat{K}$ such that $f_i'\,f=f_i$ on $U_i$ for all $i=1,\ldots,n$. (We temporarily dropped the mention of the restriction functors, for readibility purposes.)
\end{Def}

We first prove the gluing of objects without idempotent completions.
\begin{Lem}
\label{MV-obj-lem}
In a Mayer-Vietoris situation (Def.\,\ref{MV-sit-def})
$$\xymatrix{{\cat{K}}\ar[r]^-{\rho_1}\ar[d]_-{\rho_2}&{\LU{U_1}}\ar[d]^-{\rho_{12}}\\ {\LU{U_2}}\ar[r]_-{\rho_{21}}
& {\LU{U_1\cap U_2}}\,,
}$$
two objects $a_1\in\LU{U_1},\ a_2\in\LU{U_2}$ with an isomorphism $\sigma:\rho_{12} (a_1)\isoto\rho_{21}(a_2)$ in $\LU{\Ucap}$ always admit a gluing (Def.\,\ref{glue-def}).
\end{Lem}

\begin{proof}
The isomorphism $\sigma$ can be represented by a fraction $\xymatrix@C=2em{a_1& x\ar [l]_-{s}\ar[r]^-{t} & a_2}$ where $s,t$ both are $U_1\cap U_2$-isomorphisms. By Lemma \ref{factorize-lem} $s$ and $t$ factor as $s=s_1\,s_2$ and $t=t_2\,t_1$ where $s_i, t_i$ are $U_i$-isomorphisms, see the upper part of Diagram~\eqref{stu-eq}. Now complete this diagram by taking the weak push-out of $s_2$ and $t_1$:
\begin{equation}
\label{stu-eq}
\vcenter{\xymatrix@C=1.5em@R=1.5em{a_1& &x\ar[ll]_-{s}\ar[dl]^{s_2}\ar[rr]^t \ar[dr]_{t_1}& &a_2
\\
&y\ar[ul]^ {s_1}\ar[dr]_{u_1}& &z\ar[ur]_{t_2}\ar[dl]^{u_2}&
\\
& &a& & 
}}
\end{equation}
Applying Lemma~\ref{U-iso-lem}, $u_i$ is a $U_i$-isomorphism. The object $a$ is then isomorphic to $a_i$ over $U_i$ via $f_1:=s_1\circ u_1\inv$ and $f_2:=t_2\circ u_2\inv$ respectively; the relation $\sigma\,f_1=f_2$ is satisfied in $\LU{\Ucap}$ because of~\eqref{stu-eq}.
\end{proof}

\begin{Thm}[Gluing of two objects]
\label{MV-obj-thm}
In a Mayer-Vietoris situation (Def.\,\ref{MV-sit-def})
$$\xymatrix{
{\cat{K\,}}\ar[r]^-{\rho_{1}}\ar[d]_-{\rho_2}
& {\KU{U_1}}\ar[d]^-{\rho_{12}}
\\
{\KU{U_2}}\ar[r]_-{\rho_{21}}
&{\KU{\Ucap}}\,,
}
$$
given two objects $a_i\in\KU{U_i}$ for $i=1,2$ and an isomorphism $\sigma:\rho_{12} (a_1)\isoto\rho_{21}(a_2)$ in $\KU{\Ucap}$, there exists a gluing (Def.\,\ref{glue-def}), which is unique up to (possibly non-unique) isomorphism.
\end{Thm}

\begin{proof}
Obviously $\sigma\oplus T\sigma:\rho_{12}(a_1 \oplus Ta_1)\simeq\rho_{21}(a_2\oplus Ta_2)$. But $b_i:=a_i\oplus Ta_i $ is an object of $\LU{U_i}$ and since $\LU{\Ucap}\to\widetilde {\LU{\Ucap}}=\KU{\Ucap}$ is fully faithful, the morphism $\sigma\oplus T\sigma$ is an isomorphism in $\LU{\Ucap}$ as well. By Lemma~\ref{MV-obj-lem} there exists an object $b\in\cat{K}$ and isomorphisms $f_i:\rho_i(b)\to b_i$ in $\LU{U_i}$ such that $ (\sigma\oplus T\sigma)\circ\rho_{12}(f_1)=\rho_{21}(f_2)$. Consider now, for  each $i=1,2$ the morphism $\pi_i:\rho_i(b)\to \rho_i(b)$ in $\LU{U_i}$ defined by\,:
$$\xymatrix{\rho_i(b)\ar[r]^-{f_i}\ar@/_2em/[rrr]_-{\pi_i}& b_i\ar[r]^-{\smallmatrice{1&0\\0&0}}& b_i\ar[r]^-{f_i^{-1}}&\rho_i(b)\,,}$$
where $\smallmatrice{1&0\\0&0}$ on $b_i=a_i\oplus T(a_i)$ is the projection on~$a_i$\,. Now, since the diagram
$$\xymatrix@C=3.5em{\rho_{12}(b)\ar@{=}[d]\ar[r]^-{\rho_{12}(f_1)}&  \rho_{12}(b_1)\ar[d]^-{\sigma\oplus T\sigma}\ar[r]^-{\smallmatrice{1&0\\0&0}}&  \rho_{12}(b_1)\ar[d]^{\sigma\oplus T\sigma}\ar[r]^-{\rho_{12}(f_1)^{-1}}&\rho_ {12}(b)\ar@{=}[d]\\
\rho_{21}(b)\ar[r]^-{\rho_{21}(f_2)}& \rho_{21}(b_2)\ar[r]^- {\smallmatrice{1&0\\0&0}}& \rho_{21}(b_2)\ar[r]^-{\rho_{21}(f_2)^{-1}} &\rho_{21}(b)}$$
is commutative, we have that $\rho_{12}(\pi_1)=\rho_{21}(\pi_2)$. We  can now apply Mayer-Vietoris for morphisms (Thm.\,\ref{MV-mor-thm}) to $\pi_1$ and $\pi_2$ to show that there exists $\pi:b\to b$ such  that $\rho_i(\pi)=\pi_i$ for $i=1,2$. Consider now $h=\pi^2-\pi$.  Since $h=0$ in $\LU{U_i}$ it factors through objects $c_i\in\KZ {Z_i}$ as follows\,:
$$\xymatrix@C=1.5em@R=1.5em{b\ar[rr]^h\ar[rd] && b\ar[rr]^h\ar[rd] &&  b\,.\\ & c_1\ar[ru] & & c_2\ar[ru] &}$$
Since $c_1$ and $c_2$ have disjoint supports Proposition~\ref{Hom-prop}  shows that $h^2=0$. Then, by a standard trick, $p:=\pi+h-2\pi h$  satisfies $p^2=p$ and still has the property $\rho_i(p)=\pi_i$ since $\rho_i(h)=0$. Now, our category $\cat{K}$ is idempotent complete by definition, so $b$ splits as $b\simeq \Img(p)\oplus\Ker(p)$. Setting $a=\Img(p)$ gives the desired object with the property $a\simeq a_i$ in $\KU{U_i}$. Further details are left to the reader.

For uniqueness, suppose that $(a,f_1,f_2)$ and $(a',f_1',f_2')$ are two gluings. By Mayer-Vietoris for morphisms (Thm.\,\ref{MV-mor-thm}) the morphisms $f_1\inv\circ f_1'$ and $f_2\inv\circ f_2'$ glue into a morphism $a'\to a$ which must be an isomorphism (Rem.\,\ref{X-iso-rem}).
\end{proof}

\begin{Cor}[Gluing three objects]
\label{MV-3obj-cor}
Let $\Spccat{K}=U_1\cup U_2\cup U_3$ be a covering by quasi-compact open subsets. Consider three objects $a_i\in\KU{U_i}$ for $i=1,2,3$ and three isomorphisms $\sigma_{ij}:a_j\isoto a_i$ in $\KU{U_i\cap U_j}$ for $1\leq i<j\leq3$ satisfying the cocycle relation $\sigma_{12}\circ\sigma_{23}=\sigma_{13}$ in $\KU{U_1\cap U_2\cap U_3}$. Then they admit a gluing.
\end{Cor}

\begin{proof}
Note that $\Spc(\KU{U_1\cup U_2})=U_1\cup U_2=:V$ by Proposition~\ref{spc-KU-prop}. Using gluing of two objects (Thm.\,\ref{MV-obj-thm}), we can glue $a_1$ and $a_2$ into an object $b\in\KU{V}$. Using Mayer-Vietoris for morphisms (Thm.\,\ref{MV-mor-thm}) for the covering of $V\cap U_3$ by $U_{1}\cap U_3$ and $U_2\cap U_3$\,, we can construct a (non-unique) isomorphism between $b$ and $a_3$ in $\KU{V\cap U_3}$. By gluing of two objects (Thm.\,\ref{MV-obj-thm}) for the covering of $\Spccat{K}$ given by $V$ and $U_3$\,, we can now glue $b$ and $a_3$ into an object of~$\cat{K}$.
\end{proof}

\begin{Rem}
As the above proof shows, the problem that arises with three open subsets is that the isomorphism between the objects $b\in\KU{V}$ and $a_3\in\KU{U_3}$ on $V\cap U_3$ is not unique. The various choices are parametrized by $\Hom_{{\Ucap}}(Ta_1,a_2)$ but we were not able to prove that two such choices yield isomorphic gluings and we tend to believe that this is wrong in general. Nevertheless, here is a case where the gluing works for several open subsets.
\end{Rem}

\begin{Thm}[Connective gluing of several objects]
\label{conn-glue-thm}
Let $\Spccat{K}=U_1\cup \cdots\cup U_n$ be a covering by quasi-compact open subsets for $n\geq2$. Consider objects $a_i\in\KU{U_i}$ and isomorphisms $\sigma_{ji}:a_i\isoto a_j$ on $U_i\cap U_j$ satisfying the cocycle condition $\sigma_{kj}\sigma_{ji}=\sigma_{ki}$ for $1\leq i,j,k\leq n$. Assume moreover the following \emph{Connectivity Condition\,:} For any $i=2,\ldots, n$ and for any quasi-compact open $V\subset U_i$ (it is enough to take $V$ a union of intersections of $U_1,\ldots,U_n$), we suppose that\,:
\begin{equation}
\label{conn-eq}
\Hom_{V}(Ta_i, a_i)=0\,.
\end{equation}
Then there exists a gluing (Def.\,\ref{glue-def}), which is unique up to unique isomorphism.
\end{Thm}

\begin{proof}
We prove the result by induction on $n$. Let us first establish the $n=2$ case. By gluing of two objects (Thm.\,\ref{MV-obj-thm}) we only need to prove the uniqueness of the isomorphism. To see this, it suffices to prove that for two gluings $a,a'\in\cat{K}$, two (iso)morphisms $g,g':a\to a'$ which agree on $U_1$ and $U_2$ are equal. By Mayer-Vietoris for morphisms (Thm.\,\ref{MV-mor-thm}), it suffices to show that $\Hom_{{\Ucap}}(Ta,a')=0$ which follows from the Connectivity Condition~\eqref{conn-eq} and from $a\simeq a'\simeq a_2$ on $U_2$\,.

Suppose $n\geq3$ and the result known for $n-1$. Define $V=U_1\cup\ldots\cup U_{n-1}$. Since $V$ is quasi-compact, we know by Propostion~\ref{spc-KU-prop} that $\Spc(\KU{V})=V$ and we can apply the induction hypothesis to construct a gluing $b\in\KU{V}$ with isomorphisms $g_i:b\isoto a_i$ on $U_i$ for $i=1,\ldots, n-1$, such that $\sigma_{ij}g_j=g_i$ for all $1\leq i,j\leq n-1$. Consider the intersection $W:=V\cap U_n$\,, which is covered by $n-1$ quasi-compact subsets $W=(U_1\cap U_n)\,\cup\,\ldots\,\cup\,(U_{n-1}\cap U_n)$. In the category $\KU{W}$, we have two objects $b$ and $a_n$ (\ie their restrictions, of course) which are isomorphic on $U_i\cap U_n$ for $i=1,\ldots,n-1$ in a compatible way with the $\sigma_{ij}$. By uniqueness of the gluing for $n-1$, there exists a unique isomorphism $\sigma:b\isoto a_n$ on $V\cap U_n$ such that $\sigma_{in}\,\sigma=g_i$ for $i=1,\ldots,n-1$. By the $n=2$ case, we obtain the wanted gluing $a\in\cat{K}$ of $b$ and $a_n$\,, unique up to unique isomorphism. Details are left to the careful reader. Note that the uniqueness of the isomorphism $\sigma$ (at stage $n-1$) is essential for the uniqueness of the gluing~$a$ (at stage $n$).

In the above induction, we needed that if the tuple $(U_1,\ldots,U_n\,;\,a_1,\ldots,a_n)$ satisfies the Connectivity Condition~\eqref{conn-eq} for $n$, then\,:

$\smallbullet$ the tuple $(U_1\,,\ldots,U_{n-1}\,;\,a_1,\ldots,a_{n-1})$ satifies~\eqref{conn-eq} for $n-1$,

$\smallbullet$ the tuple $(U_1\cap U_n\,,\ldots,U_{n-1}\cap U_n\,;\,a_1,\ldots,a_{n-1})$ satifies~\eqref{conn-eq} for $n-1$,

$\smallbullet$ the 4-uple $(U_1\cup \ldots \cup U_{n-1}\,,U_n\,;\,b,a_n)$ satifies~\eqref{conn-eq} for $n=2$, for any object $b$.

\noindent These are easy to check. The last one comes from the assumption $i>1$ in~\eqref{conn-eq}.
\end{proof}


\bigbreak
\section{Picard groups}
\medbreak


%
\begin{Def}
An object $a\in\cat{K}$ is called \emph{invertible} if there  exists an object $b$ such that $a\otimes b\simeq\unit$. By adjunction, see Def.\,\ref{spc-def}, an object $a\in\cat{K}$ is invertible if and only if the evaluation map $\eta:Da\otimes a\to\unit$ is an isomorphism.
\end{Def}

\begin{Lem}
\label{local-inv}
An object $a$ in $\cat{K}$ is invertible if and only if it is invertible in $\cat{K}/\cat{P}$ for all $\cat{P}\in\Spccat{K}$.
\end{Lem}

\begin{proof}
This is clear by Proposition~\ref{iso-local-prop}.
\end{proof}

\begin{Def}
Define $\Piccat{K}$ to be the set of isomorphism classes of  invertible objects in $\cat{K}$. The tensor product $\otimes:\cat{K} \times\cat{K}\to\cat{K}$ makes $\Piccat{K}$ into an abelian group  with unit the class of~$\unit$. 
\end{Def}

\begin{Prop}
\label{pic=pic-prop}
Let $X$ be a scheme and consider $\Dperf(X)$ its derived  category of perfect complexes. Then there is a split short exact sequence of abelian groups
$$0\to\Piccat{K}\to\Pic(\Dperf(X))\to C(X;\bbZ)\to 0$$
where $C(X;\bbZ)$ stands for the group of locally constant functions from $X$ to $\bbZ$.
\end{Prop}

\begin{proof}
We first describe $\Pic(\Dperf(X))$ where $X=\Spec(R)$ is the  spectrum of a local ring $(R,\mm)$. In this case, any object of $ \Dperf(R)$ is isomorphic to a so-called \emph{minimal} complex of the  form
$$C=\xymatrix{\cdots\ar[r]& R^{\ell_i}\ar[r]^-{d_i} & R^{\ell_{i-1}} \ar[r]& \cdots}$$
where, for all $i$, the differential $d_i$ is a matrix with  coefficients in $\mm$. If $C$ is invertible in $\Dperf(R)$ so is $\bar {C}$, its image under the functor $\Dperf(R)\to\Dperf(R/\mm)$. But all the differentials of $\bar{C}$ are 0 and the relation $C\otimes D \simeq R$, for some complex $D$, now shows that the complex $C$ must be $R$ concentrated in some degree, \ie there exists $n_0=n_0 (C)$ such that $\ell_i=1$ if $i=n_0$ and $\ell_i=0$ otherwise.

For a global~$X$, the map $\Pic(\Dperf(X))\to C(X;\bbZ)$ is now easily defined: for an  invertible complex $C\in\Dperf(X)$ and for $x\in X$ denote by $C_x$  its image in $\Dperf(\calO_{X,x})$. The function $f_C:X\to\bbZ$ is  then defined by $x\mapsto n_0(C_x)$. The rest of the proof is straightforward\,: a perfect complex which is locally trivial is quasi-isomorphic to its homology in degree zero and the latter must be a line bundle.
\end{proof}

\begin{Def}
\label{Gm-def}
Define $\Gmcat{K}=\Hom_{\cat{K}}(\unit,\unit)^{\times}$ to be the group of invertible elements of the  (commutative) ring $\Hom_{\cat{K}}(\unit,\unit)$.
\end{Def}

\begin{Thm}
\label{MV-pic-thm}
In a Mayer-Vietoris situation (Def.\,\ref{MV-sit-def}), there is an exact sequence of abelian  groups
\vskip-\baselineskip
\vskip-\baselineskip
$$\xymatrix@C=2em@R=1em{&&\kern8em\cdots\ar[r]&\Homcat{K(\Ucap)}(T\unit,\unit)\ar[r]^-{1+ \partial}&
\\
\ar[r]^-{1+\partial}&{\Gmcat{K}}\ar[r]&{\Gm}(\KU{U_1})\oplus{\Gm}(\KU{U_2})\ar[r]
&{\Gm}(\KU{\Ucap})\ar[r]^-{\delta}&
\\
\ar[r]^-{\delta}&\Piccat{K}\ar[r]&\Pic(\KU{U_1})\oplus\Pic(\KU{U_2})\ar[r]& \Pic(\KU{\Ucap})\,.}$$
The homomorphism $\partial$ is as in Construction~\ref{partial-cons} and the unlabelled homomorphisms are the restrictions and the (multiplicative) differences thereof.

The homomorphism $\delta:\Gm(\KU{\Ucap})\to \Piccat{K}$ is defined by gluing two copies of $\unit$ by means of Theorem~\ref{MV-obj-thm}. Explicitly, it can be described as follows\,: Let $\sigma=ts\inv:\ \xymatrix@C=2em{\unit& x\ar[l]_-{s}\ar[r]^-{t} & \unit}$ in $\Gm(\KU{\Ucap})$, where $s$ and  $t$ are $\Ucap$-isomorphisms\,; by Lemma \ref{factorize-lem}, there exist factorizations $s=s_1\,s_2$ and $t=t_2\,t_1$ where $s_i$ and $t_i$ are $U_i$-isomorphisms\,; then $\delta(\sigma)$ is defined as the isomorphism class of the weak push-out $p\in\cat{K}$ of $x_1$ and $x_2$ over $x$
\begin{equation}
\label{delta-eq}
\vcenter{\xymatrix@C=1.5em@R=1.5em{\unit& &x\ar[ll]_-{s}\ar[dl]^{s_2}\ar[rr] ^t\ar[dr]_{t_1}& &\unit
\\
&x_1\ar[ul] ^{s_1}\ar[dr]_{u_1}& &x_2\ar[ur]_{t_2}\ar[dl]^{u_2}&
\\
& &p\,.\kern-0.4em
}
}\end{equation}
\end{Thm}

\begin{proof}
First note that the homomorphism $1+\partial:\Hom_{\Ucap}(T\unit,\unit)\to\Homcat{K}(\unit,\unit)$ lands in $\Gmcat{K}$. Indeed for any $g\in\Hom_{\Ucap}(T\unit,\unit)$ one has $\partial(g)\circ\partial(g)=0$, since $\partial(g):\unit\to\unit$ is zero in $\KU{U_i}$ and hence factors via $\KZ{Z_i}$ for $i=1,2$ and since $\Homcat{K}(\KZ{Z_1},\KZ{Z_2})=0$ by Prop.\,\ref{Hom-prop}. So, $1+\partial(g)$ is invertible with inverse $1-\partial(g)$.

The connecting homomorphism $\delta:\Gm(\KU{\Ucap})\to\Piccat{K}$ produces an object $p\in\cat{K}$, see Diagram~\eqref{delta-eq}, which is isomorphic to $1$ on $U_1$ via $s_1\,u_1\inv$ and on $U_2$ via $t_2\,u_2\inv$, in a compatible way with $\sigma$ on~$\Ucap$. The object~$p$ is then the gluing of two copies of $1$ along the isomorphism $\sigma$ on $\Ucap$. Such a gluing is unique up to isomorphism by gluing of two objects (Thm.\,\ref{MV-obj-thm}), and this gluing only depends on the map~$\sigma$, and not on the various choices ($s,t,s_1,s_2,t_1,t_2,p$). So, the map $\delta$ is well-defined and we now check that it is a group homomorphism. Take $\xymatrix@C=2em{\unit& x\ar[l]_-{s}\ar[r]^-{t} & \unit}$ and $ \xymatrix@C=2em{\unit& x'\ar[l]_-{s'}\ar[r]^-{t'} & \unit}$ two units in $\Gm(\KU{\Ucap})$. With the same notations as above, factor these morphisms and perform the corresponding weak push-outs as follows
$$
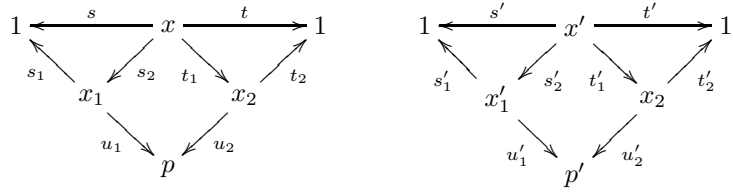
\xymatrix@C=1.5em@R=1.5em{\unit& &x\ar[ll]_-{s}\ar[dl]^{s_2}\ar[rr] ^t\ar[dr]_{t_1}& &\unit\\
                                                          &x_1\ar[ul] ^{s_1}\ar[dr]_{u_1}& &x_2\ar[ur]_{t_2}\ar[dl]^{u_2}& \\
                    & &p& & }{}\hskip1cm \xymatrix@C=1.4em@R=1.3em {\unit& &x'\ar[ll]_-{s'}\ar[dl]^{s'_2}\ar[rr]^{t'}\ar[dr]_{t'_1}& & \unit\\
                                                          &x'_1\ar [ul]^{s'_1}\ar[dr]_{u'_1}& &x_2\ar[ur]_{t'_2}\ar[dl]^{u'_2}& \\
                    & &p'& & }{}$$
In a symmetric monoidal category, the composition of two morphisms between the unit object is also given by the tensor product. Hence we tensor together the two above diagrams to obtain the following one\,:
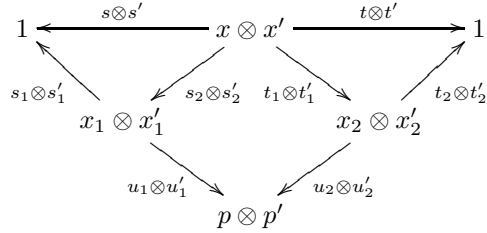
$$\xymatrix@C=1.3em@R=2em{\unit& &x\otimes x'\ar[ll]_-{s\otimes s'}\ar [dl]^{s_2\otimes s'_2}\ar[rr]^-{t\otimes t'}\ar[dr]_{t_1\otimes t'_1} & &\unit\\
                                                          &x_1 \otimes x'_1\ar[ul]^{s_1\otimes s'_1}\ar[dr]_{u_1\otimes u'_1}& &x_2 \otimes x'_2\ar[ur]_{t_2\otimes t'_2}\ar[dl]^{u_2\otimes u'_2}& \\
                    & &p\otimes p'& & }{}$$
By Lemma~\ref{bicar-lem}, the above middle square is weakly bicartesian  as well. Hence, $p\otimes p'=\delta(\sigma\otimes \sigma')=\delta(\sigma\,\sigma')$. This shows that $\delta$ is an group homomorphism.

(Recall the restriction functors $\rho_i$ and $\rho_{ij}$ from Definition~\ref{MV-sit-def}.) It is straightforward from the above definition of $\delta$ that $\rho_i\circ\delta=1$ for $i=1,2$. To see that $\delta\circ\rho_{12}=1$, for instance, one can assume that $s_2=\id$ and $t_2=\id_{\unit}$ in~\eqref{delta-eq}, in which case $u_2$ must also be an isomorphism, \ie $p\simeq x_2=\unit$. The other compositions are clearly 1. The exactness of the left-hand side of the sequence up to $\Gmcat{K}$ follows from Mayer-Vietoris for morphisms (Thm.\,\ref{MV-mor-thm}) applied at $a=b=\unit$. It remains to check the exactness of the sequence at four spots.
\medbreak
\noindent \textit{Exactness at $\Gm(\KU{U_1})\oplus\Gm(\KU{U_2})$}: This is again immediate from Mayer-Vietoris for morphisms (Thm.\,\ref{MV-mor-thm}) recalling that a local isomorphism is an isomorphism (Prop.\,\ref{iso-local-prop}).
\medbreak
\noindent \textit{Exactness at $\Gm(\KU{\Ucap})$}: Let $\sigma=\big(\xymatrix@C=2em{\unit& x\ar[l]_-{s}\ar[r]^-{t} & \unit}\big)$ in $\Gm(\KU{\Ucap})$ be such that $\delta(\sigma)\simeq\unit$ in $\cat{K}$. This means that one can find a diagram of the form
$$\xymatrix@C=1.5em@R=1.5em{\unit& &x\ar[ll]_-{s}\ar[dl]^{s_2}\ar[rr] ^t\ar[dr]_{t_1}& &\unit\\
                                                          &x_1\ar[ul] ^{s_1}\ar[dr]_{u_1}& &x_2\ar[ur]_{t_2}\ar[dl]^{u_2}& \\
                    & &\unit& & }{}$$
see~\eqref{delta-eq}. One then sees two morphisms, namely $\sigma_1=u_1\,s_1\inv\in\Gm(\KU{U_1})$ and $\sigma_2=u_2\,t_2\inv\in\Gm(\KU{U_2})$, such that $\sigma_2^{-1}\circ \sigma_1=\sigma$ in $\Gm(\KU{\Ucap})$.
\medbreak
\noindent \textit{Exactness at $\Piccat{K}$}: Let $p$ be an invertible object in $\cat{K}$ such that $\rho_i(p)\simeq\unit$ for $i=1,2$. Thus there exist $U_i$-isomorphisms $\xymatrix@C=2em{\unit&  y_i\ar[l]_-{s_i}\ar[r]^-{t_i} & p}$. Performing the weak pull-back of $s_1$ and $s_2$ one obtains the diagram
$$
\xymatrix@C=1.5em@R=1.5em{\unit& &y\ar[dl]_{u_2}\ar[dr]^{u_1}&&\unit
\\
&y_1\ar[ul]^{t_1}\ar[dr]_{s_1}& &y_2\ar[ur]_{t_2}\ar[dl]^{s_2}&
\\
& &p
}
$$
and this defines $\big(\xymatrix@C=2.5em{\unit& y\ar[l]_-{t_1\, u_2}\ar[r] ^-{t_2\, u_1} & \unit}\big)\in\Gm(\KU{\Ucap})$. The image of this morphism under $\delta$ is clearly isomorphic to $p$ by construction, see~\eqref{delta-eq}, the middle square in the above diagram being weakly bicartesian.
\medbreak
\noindent \textit{Exactness at $\Pic(\KU{U_1})\oplus\Pic(\KU{U_2})$}: This follows from gluing of two objects (Thm.\,\ref{MV-obj-thm}) and from invertibility being a local property (see Lemma~\ref{local-inv}).
\end{proof}

\begin{Thm}
\label{H1-thm}
Let $\cat{K}$ be an \icsctenscat. Suppose that $\Hom_{{U}}(T\unit,\unit)=0$ for every quasi-compact open subsets~$U\subset\Spccat{K}$. Then there exists a unique sheaf $\Gm$ on $\Spccat{K}$, such that $\Gm(U)=\Gm(\KU{U})$ when $U\subset\Spccat{K}$ is quasi-compact open. Moreover, there exists an injective homomorphism from the first \v Cech cohomology of $\Spccat{K}$ with coefficients in~$\Gm$ into the Picard group of~$\cat{K}$
$$
\alpha\,:\cH^1(\Spccat{K},\Gm)\hookrightarrow \Piccat{K}
$$
which sends a $\Gm$-cocycle~$\sigma$ to the unique gluing of copies of $\unit$ along the isomorphisms over the pairwise intersections given by~$\sigma$, as described in Theorem~\ref{conn-glue-thm}.
\end{Thm}

\begin{proof}
We first prove by induction on~$n$ the following

\textit{Claim}\,: Given a covering of a quasi-compact subset $V\subset\Spccat{K}$ by $n\geq 2$ quasi-compact open subsets, $V=U_1\cup\ldots\cup U_n$, and given morphisms $f_{i}:\unit\to \unit$ over $U_i$ such that $f_i=f_j$ over $U_i\cap U_j$ for $1\leq i,j\leq n$, there exists a unique morphism $f:\unit \to \unit$ over~$V$ such that $f=f_i$ over $U_i$\,.

Replacing $\cat{K}$ by $\KU{V}$, we can assume that $V=\Spccat{K}$. Now, for $n=2$, this is Mayer-Vietoris for morphisms (Thm.\,\ref{MV-mor-thm}). Note that unicity follows from $\Hom_{{\Ucap}}(T\unit,\unit)=0$. The induction on~$n$ is then easy\,: To construct $f$, glue the $n-1$ first morphisms $f_i$ into a morphism $g:\unit\to\unit$ on $U_1\cup\ldots\cup U_{n-1}$ and show that it agrees with $f_n$ on the intersection with $U_n$ -- this uses unicity for $n-1$\,; then apply the $n=2$ case to glue $g$ and $f_n$ into a global $f$. To prove unicity of~$f$, proceed similarly, using unicity for $n-1$ and for $n=2$ again. Hence the Claim.

\smallbreak

Then the existence of the sheaf~$\Gm$ is immediate from the claim and from the fact that quasi-compact open subsets form a basis of the topology of $\Spccat{K}$ (\cite{b05} Rem.\,2.7 and Prop.\,2.14). For the same reason and because of quasi-compacity of $\Spccat{K}$, to define the homomorphism $\alpha$, it suffices to consider $\Gm$-cocycles over finite coverings of $\Spccat{K}$ by quasi-compact open subsets. In this situation, the gluing is guaranteed by Theorem~\ref{conn-glue-thm}. Hence $\alpha$ is well-defined.

\smallbreak

Finally, injectivity of $\alpha$ is easy. Indeed, given a $\Gm$-cocyle $\sigma$ over a covering $\Spccat{K}=U_1\cup\ldots\cup U_n$ with every $U_i$ quasi-compact open, the gluing $a\in\Piccat{K}$ comes with isomorphisms $f_i:a\isoto \unit$ on each $U_i$, compatible with the $\sigma(U_i\cap U_j)$ as usual. Now, if $a=\unit$, the latter compatibilty means that the \v Cech boundary of the 0-cochain defined by the $f_i\in\Gm(U_i)$ is nothing but $\sigma$, that is, $\sigma=0$ in $\cH^1(\Spccat{K},\Gm)$.
\end{proof}

\begin{Rems}
\label{pic-rem}
\begin{enumerate}
\item
Note that the condition $\Hom_{{U}}(T\unit,\unit)=0$ does not hold in general, for instance in Modular Representation Theory, \ie for $\cat{K}=kG\stab$. For instance, for $k=\bbF_2$ and $G=\bbZ/2$, we even have $T\unit\simeq\unit$\,!
\smallbreak
\item
When the condition $\Hom_{{U}}(T\unit,\unit)=0$ holds for every quasi-compact open $U\subset\Spccat{K}$ and when $\Spccat{K}$ happens to be a scheme, Theorem~\ref{H1-thm} gives an injective homomorphism $\Pic(\Spccat{K})\hookrightarrow\Piccat{K}$. In the case of $\cat{K}=\Dperf(X)$ for $X$ a scheme, this homomorphism is the one of Proposition~\ref{pic=pic-prop}.
\smallbreak
\item
The previous remark clearly leads us to glue copies of $T^r(\unit)$ for any $r\in\bbZ$, or even for any locally constant function $r\in C(\Spccat{K},\bbZ)$. This induces a homomorphism\,:
$$
\cH^1(\Spccat{K},\Gm)\times C(\Spccat{K},\bbZ)\too \Piccat{K}
$$
which we do not know to be injective in general.
\end{enumerate}
\end{Rems}


\bigbreak
\section{Excision}
\medbreak


For later use, we state the following result in greater generality than in the Introduction. See Remark~\ref{more-rem}. The reader can as well consider the case of $A$ and $B$ reduced to a singleton, \ie $U$ open and $Y$ closed.

\begin{Thm}[Excision]
\label{exc-thm}
Let $\cat{K}$ be an \icsctenscat\ and let $Y\subset U\subset \Spccat{K}$. Assume that $Y=\cup_{\alpha\in A} Y_\alpha$ with every $Y_\alpha$ closed with quasi-compact complement and assume that $U=\cap_{\beta\in B}U_\beta$ with every $U_\beta$ open and quasi-compact. Then the restriction functor $\rho:\cat{K}\to\KU{U}$ induces an equivalence between the respective subcategories of objects supported on~$Y$\,:
$$
\KZ{Y}\isotoo\KU{U}_{Y}\,.
$$
\end{Thm}

\begin{proof}
Remark first of all that $\Spc(\KU{U})\cong U$ by Prop.\,\ref{spc-KU-prop}, whose proof generalizes \textsl{verbatim} to this situation.

Let us see that the functor $\rho:\KZ{Y}\to \KU{U}_Y$ is full. Given $a,b\in\KZ{Y}$ and a fraction $\xymatrix@C=2em{a& x\ar[l]_-{s}\ar[r]^-{f} & b}$ with $s$ a $U$-isomorphism, we have $\supp(\cone(s))\cap \supp(Ta)\subset \supp(\cone(s))\cap U=\varnothing$, so $\Homcat{K}(a,\cone(s))=0$ by Prop.\,\ref{Hom-prop}. So, the morphism~$s$ is a split epimorphism, say $s\circ u=\id_a$ for a morphism $u:a\to x$. Amplifying the fraction $f\,s\inv$ by $u$ shows that this morphism $f\,s\inv$ is equal to (the restriction of) the morphism $fu:a\to b$.

Let us see that the functor $\rho:\KZ{Y}\to \KU{U}_Y$ is faithful. Let $f:a\to b$ be a morphism in $\KZ{Y}$ such that $\rho(f)=0$, that is, there exists a $U$-isomorphism $s:x\to a$ such that $fs=0$. As above, $s$ must be a split epimorphism, hence $f=0$.

Let us see that the functor $\rho:\KZ{Y}\to \KU{U}_Y$ is essentially surjective. Let $b\in\KU{U}_Y$. There is an object $a\in\cat{K}$ such that $\rho(a)=b\oplus T(b)$. We have $\supp(a)\cap U=\supp_U(b)\cup\supp_U(Tb)=\supp_U(b)\subset Y$. So, if we call $Z=\Spccat{K}\oursetminus U$ the complement of $U$, we have proved that $\supp(a)\subset Y\cup Z$. Since $Y\cap Z=\varnothing$, we know by Theorem~\ref{Balmson} that $a\simeq c\oplus d$ with $\supp(c)\subset Y$ and $\supp(d)\subset Z$. But then $\rho(a)=\rho(c)$ and we have found an object $c\in\KZ{Y}$ such that $\rho(c)=b\oplus T(b)$. Now, consider the idempotent of $b\oplus T(b)$ corresponding to the projection on~$b$. Since $\rho$ is fully faithful, there exists a corresponding idempotent on the object $c$, which then decomposes accordingly, one factor going to $b$, as was to be shown.
\end{proof}

\begin{Rem}
If needed, the reader can formalize the following assertion\,: Given a point $\cat{P}\in\Spccat{K}$, the ``local'' category $\cat{K}/\cat{P}$, or rather its idempotent completion, is equivalent to the colimit of the categories $\KU{U}$, over the quasi-compact open neighborhoods $U\ni\cat{P}$.
\end{Rem}


\bigbreak
\noindent\textbf{Acknowledgments\,:}
We thank Jacques Th\'evenaz for a counter-example, which prevented us from losing too much time on some wild dream of ours. The second author thanks the FIM of the Mathematics Department at the ETH Zurich for impartial partial support.


\bibliographystyle{abbrv}
\bibliography{TG-articles}


\end{document}